\documentclass[12pt]{amsart}
\usepackage{latexsym,amssymb,amsmath,amsthm}
\usepackage{comment}
\usepackage{float}
\usepackage{graphicx}
\newtheorem{thm}{Theorem}

\theoremstyle{definition}

     \def\Det{\rm{Det}} \def\det{{\rm det}} 
\def\tr{\rm{tr}}

\title{Counting rooted forests in a network}
\author{Oliver Knill}
\date{July 13, 2013, Update July 18, 2013}
\address{
        Department of Mathematics \\
        Harvard University \\
        Cambridge, MA, 02138
        }

\begin{document}
\maketitle
\begin{abstract}
If $F,G$ are two $n \times m$ matrices, then $\det(1+x F^T G) = \sum_P x^{|P|} \det(F_P) \det(G_P)$ 
where the sum is over all minors \cite{cauchybinet}. An application is a new proof of the
Chebotarev-Shamis forest theorem telling that $\det(1+L)$ is the number of rooted 
spanning forests in a finite simple graph $G$ with Laplacian $L$. We can generalize this and show
that $\det(1+k L)$ is the number of rooted edge-k-colored spanning forests.  If a forest with an even number
of edges is called even, then $\det(1-L)$ is the difference between even and odd rooted spanning 
forests in $G$. 
\end{abstract} 

\section{The forest theorem}

A social network describing friendship relations is mathematically described by 
a finite simple graph.  Assume that everybody can chose among their friends a candidate
for ``president" or decide not to vote. How many possibilities are there to do so, if 
cyclic nominations are discarded? The answer is given explicitly as the product of $1+\lambda_j$,
where $\lambda_j$ are the eigenvalues of the combinatorial Laplacian $L$ of $G$. More generally, 
if votes can come in $k$ categories, then the number voting situation is the product of
$1+k \lambda_j$. 
We can interpret the result as counting rooted spanning forests in finite simple graphs, which is 
a theorem of Chebotarev-Shamis. 
In a generalized setup, the edges can have $k$ colors and get a formula for these rooted spanning 
forests.  While counting subtrees in a graph is difficult 
\cite{Jerrum,GoldbergJerrum} in Valiants complexity class $\#P$, Chebotarev-Shamis show that this
is different if the trees are rooted. 
The forest counting result belongs to spectral graph theory \cite{Biggs,Chung97,CDS,Beineke,Knauer} 
or enumerative combinatorics \cite{HararyPalmer,Godsil}. 
Other results relating the spectrum of $L$ with combinatorial properties is
Kirchhoff's matrix tree theorem which expresses the number of spanning trees in a 
connected graph of $n$ nodes as the pseudo determinant $\Det(L)/n$ or the Google determinant
$\det(E+L)$ with $E_{ij}=1/n^2$. 
counting the number of rooted spanning trees in $G$, a measure for complexity of the graph \cite{Biggs}.
Of course, the number of people in the network 
is $\tr(L^0)$ and by the handshaking lemma of Euler, the number of friendships is $\tr(L)/2$. 
An other example of spectral-combinatorial type is that the 
largest eigenvalue $\lambda_1$ of $L$ gives the upper bound $\lambda_1-1$ for the maximal 
number of friends which can occur. A rather general relation between the sorted eigenvalues 
$\lambda_1 \geq \lambda_2 \dots $ and degrees $d_1 \geq d_2 \geq \dots$ is 
$\sum_{i=1}^k \lambda_i \geq \sum_{i=1}^k d_i$ (\cite{CRS} Theorem 7.1.3).
While working on matrix tree theorems for the Dirac operator \cite{DiracKnill}, we have found a 
generalization of the classical Cauchy-Binet theorem in linear algebra \cite{cauchybinet}.
It tells that if two $m \times n$ matrices $F,G$ of the same size are given, then 
two polynomials agree: one is the characteristic polynomial of the $n \times n$ 
matrix $F^TG$ and the other is polynomial containing product of all possible minors of $F$ or $G$:
\begin{equation}
 p(x) = \det(F^T G-x) = \sum_{k=0}^n (-x)^{n-k} \sum_{|P|=k} \det(F_P) \det(G_P) \; ,
 \label{1}
\end{equation}
Here $\det(H_P)$ is the minor in $H$ masked by a square pattern $P = I \times J$ and $x= 1_n x$ is
the diagonal matrix with entry $x$ and for $k=0$, the summand with the empty pattern
$P=\emptyset$ is understood to be $1$. 
The classical Cauchy-Binet formula is the special case, when $x=0$ and $F^TG$ is invertible. 
The proof of formula~(\ref{1}) is given in \cite{cauchybinet} using
exterior calculus. While multilinear algebra proofs of Cauchy-Binet have entered
textbooks \cite{HuppertWillems,Shafarevich}, the identity (\ref{1}) appears to be new.
We have not yet stated in \cite{cauchybinet} that for $x=-1$, we get the remarkably general 
but new formula for classical determinants
\begin{equation}
 \det(1+F^T G) = \sum_P \det(F_P) \det(G_P)  \; , 
\label{2} 
\end{equation}
where the sum is over all minors and which is true for all matrices $F,G$ of 
the same size and where the right hand side is $1$ if $P$ is empty. 
If $F,G$ are column vectors, then the identity tells $1+\langle F,G \rangle = 1+\sum_i F_i G_i$
so formula~(\ref{2}) generalizes the dot product. 
For square matrices $A$, it implies the Pythagorean identity
\begin{equation}
 \det(1+A^T A) = \sum_P \det^2(A_P)  \; ,
\label{3} 
\end{equation}
where again on the right hand side the sum is over all minors. Even this special case seems have 
been unnoticed so far. In the graph case, where $L=C^T C={\rm div} \circ {\rm grad}$ 
for the incidence matrix $C$ = ``gradient", formula ~(\ref{3})  implies 
for $F=C,G=C^T$ the relation $\det(1+L) = \sum_P \det^2(C_P)$ for the Laplacian $L$. 
Poincar\'e has shown in 1901 that $\det^2(C_P) \in \{0,1 \;\}$. Actually, it is $1$ if and only if $P$ 
belongs to a subchain of the graph obtained by choosing the same number of edges and vertices in such a way
that every edge connects with exactly one vertex and so that we do not form loops. These are rooted
forests, collections of rooted trees. Trees with one single vertex are seeds which when include lead to
rooted spanning forests.  From formula (\ref{3}) follows a theorem of Chebotarev and Shamis
(which we were not aware of when first posting the result):

\begin{thm}[Chebotarev-Shamis Forest Theorem]
For a finite simple graph $G$ with Laplacian $L$, the integer $\det(1+L)$ is the number 
of rooted spanning forests contained in $G$. 
\end{thm}

With the more general formula
\begin{equation}
 \det(1+k A^T A) = \sum_P k^{|P|} \det^2(A_P)  \; ,
\label{4} 
\end{equation}
the forest theorem a bit further and get a more general result which counts 
forests in which branches are colored. Lets call a graph $k$-colored, if its
edges can have $k$ colors: 

\begin{thm}[Forest Coloring Theorem]
For a finite simple graph $G$ with Laplacian $L$, the integer $\det(1+kL)$ is the number
of rooted $k$-colored spanning forests contained in $G$.
\end{thm}

We can also look at $k=-1$, in which case we count forests with odd number of trees
with a negative sign. 
Lets call a graph ``even" if it has an even number of edges,
and ``odd" if it has an odd number of edges. 

\begin{thm}[Super Forest Coloring Theorem]
For a finite simple graph $G$ with Laplacian $L$, the integer
$\det(1-kL)$ is the number of $k$-colored rooted even spanning forests 
minus the number of $k$-colored rooted odd spanning forests in $G$.
\end{thm}

\section{Remarks}

The integer $\det(1+L)$ is also the number of simple outdegree=1 acyclic digraphs contained in $G$. 
It is the number of voting patterns excluding mutual and cyclic votes in a finite network
or the number of molecule formations with $n$ atoms, where each molecule is of the
form $4C_n H_{2n+1} X$ with $X$ representing any radical different from hydrogen \cite{Otter}
representing the root. In the voting picture, the root is a person in the tree which is voted on, 
but does not vote. Rooted trees are pivotal in computer science, because directories are
rooted trees. A collection of virtual machines can be seen as a rooted forest.
The colored matrix forest theorem can have an interpretation also in that there are $k$
issues to vote for and that the root person who does not vote can chose the issue. 
As mentioned the {\bf forest theorem} is due to Chebotarev-Shamis 
\cite{ChebotarevShamis2,ChebotarevShamis1}.
There are enumeration results \cite{ErdoesForests}
and generating functions \cite{Lampe} motivated by the Jacobian conjecture or
tree packing results \cite{Katoh}. 
Chung-Zhao \cite{ChungZhao} have a matrix forest theorem for the normalized Laplacian which
involves a double sum with weights involving the degrees of the vertices.
The double sum is over the number of trees as well as the roots.  \\

The integer valued function $f(G) = \det(1+L)$ on the category of finite simple graphs is positive since
$L$ has nonnegative spectrum. The combinatorial description immediately implies that $f$ is monotone
if $k>0$: if $H$ is a subgraph of $G$, then $f(H) \leq f(G)$. 
Knowing all the numbers $\det(L+k)$ does not characterize a graph even among connected graphs 
as many classes of isospectral graphs are known. 
Since we know $f$ explicitly on complete graphs and graphs without edges,
it follows that $1 \leq f(G) \leq (1+k n)^{n-1}$ if $n$ is the order of the graph and
$k$ colors are used.  Finally, lets mention that if $G$ is the disjoint union of two 
graphs $G_1$ and $G_2$, then $f(G) = f(G_1) f(G_2)$ for any $k$. This is both clear combinatorially
as well as algebraically, because the graph $G$ has as the eigenvalues the union of the eigenvalues
of $G_1$ and $G_2$. 

\section{Examples}

Extreme cases are the zero dimensional graph with $n$ vertices and no edges as well as the 
complete graph $K_n$. In the zero dimensional case, all trees are points and there is one
possibility of a maximal spanning forest. Already in the $K_n$ case, a direct computation
using partitions and applying the Cayley formula $n_j^{n_j-2}$ for the number of spanning 
trees in each subset of cardinality $n_j$ is not easy. The total sum $f(K_n) = (n+1)^{n-1}$
looks like the Cayley formula but this is a coincidence because all nonzero eigenvalues 
of $L$ are $n$ so that $\Det(L)/n$ and $\det(L+1)$ look similar. The following examples
are for $k=1$ where we count the number of spanning forests. \\

{\bf 1)} Zero dimensional: no connections $f(G) = 1$. \\
{\bf 2)} Complete $f(K_n) = (n+1)^{n-1}$ like $f(K_2)=3,f(K_3)=16$. \\
{\bf 3)} Star graphs: $f(S_n) = (n+1) 2^{n-2}$ like $f(S_2)=3, f(S_3)=8$. \\
{\bf 4)} Cyclic: $f(C_n) = \prod_k (1+4 \sin^2(\frac{\pi k}{n})$ like $f(C_3)=16, f(C_4)=45$. \\
{\bf 5)} Line graph: $f(L_n)$ =bisected Fibonacci: $f(L_2)=3,f(L_3)=8$ ...\\
{\bf 6)} Wheel graph:  $f(3)=125,f(4)=576,f(5)=2527,f(6)=10800$. \\ 
{\bf 7)} Bipartite: $f(1)=3,f(2)=45,f(3)=1792,f(4)=140625$. \\
{\bf 8)} Platonic: $f({\rm T,O,H,D,I})=(125, 6125, 23625, 500697337, 107307307008)$. \\
{\bf 9)} Molecules: $f({\rm caffeine}) =7604245376$, $f({\rm guanine}) = 0$ \\
{\bf 10)} Random Erdoes-Renyi: $f$ appears asymptotically normal on $E(n,p)$. \\

Actually, \cite{Chebotarev2008} analyzed the line graph case and confirmed the 
Fibonnacci connection. We had only noticed this experimentally.  \\
According to \cite{Callan}, Cayley knew also the number of rooted forests in 
$K_n$ as $(n+1)^{n-1}$. \\

In the following examples, where $k=-1$ and where we count the difference between
the odd and even rooted forests: \\

{\bf 1)} Complete $f(K_n) = (1-n)^{n-1}$ like $f(K_2)=-1,f(K_3)=4$. \\
{\bf 2)} Cyclic: $f(C_n) = \prod_{k=1}^{n-1}(1-4 \sin^2(\frac{\pi k}{n}))$ is 6 periodic.  \\
{\bf 3)} Star and line graphs: $f(S_n) = f(L_n) =  0$. \\
{\bf 4)} Platonic: $f({\rm T,O,H,D,I})=(-135, 4096, -4159375, -675, -27)$. \\
{\bf 5)} Erdoes-Renyi: $f$ appears asymptotically normal on $E(n,p)$. \\

Remarkable is that for cyclic graphs, also the sequence $\det(L(C_n)-1)$ is cyclic. We have
$f(n) = \prod_{k=1}^{n-1} (4 \sin(\pi k/n)^2-1)$ and this is always $6$ periodic in $n$:
$f(1)=1,f(2)=3,f(3)=4,f(4)=3,f(5)=1,f(6)=0,f(7)=1, \dots$. 
For  star and line graphs $\det(L(C_n)-1)=0$ because they have an eigenvalue $1$.
Interesting that for Platonic solids, the cube has the most 
deviation from symmetry, the octahedron is the only positive 
and that the icosahedron has maximal symmetry between odd and even forests. \\

In the following examples, the seeds, trees with one vertex only, are not marked.
The figures illustrate also the voting picture (our first interpretation). 
For these illustrations, we also assume that $k=1$, we do not color the forests. \\

{\bf Example 1.}
The triangle with Laplacian
$L=\left[
\begin{array}{ccc}
 2 & -1 & -1 \\
 -1 & 2 & -1 \\
 -1 & -1 & 2 \\
\end{array} \right]$ has the eigenvalues of $L$ are $0,3,3$ so that $f(G)=16$. 

\begin{figure}
\scalebox{0.26}{\includegraphics{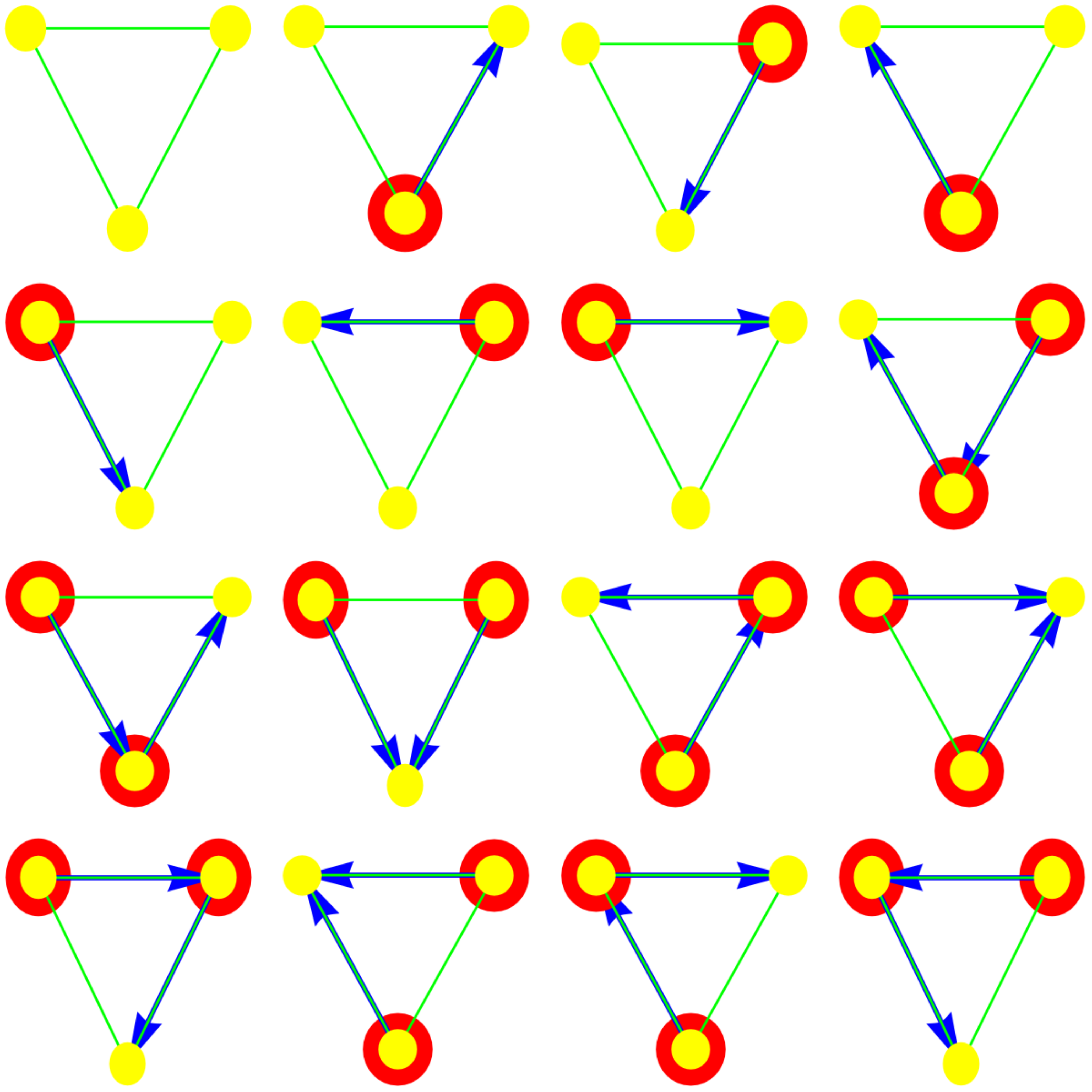}}
\end{figure}

{\bf Example 2.}
A graph $G$ with Laplacian
$L=\left[
\begin{array}{cccc}
 3 & -1 & -1 & -1 \\
 -1 & 2 & -1 & 0 \\
 -1 & -1 & 2 & 0 \\
 -1 & 0 & 0 & 1 \\
\end{array} \right]$ is called $Z_1$ \cite{HHM}. 
The eigenvalues of $L$ are $0,1,3,4$ so that $f(G) =40$. 
The last $12 = \Det(L)$ are rooted spanning trees. 
We see here already examples with two disjoint trees. 

\begin{figure}
\scalebox{0.32}{\includegraphics{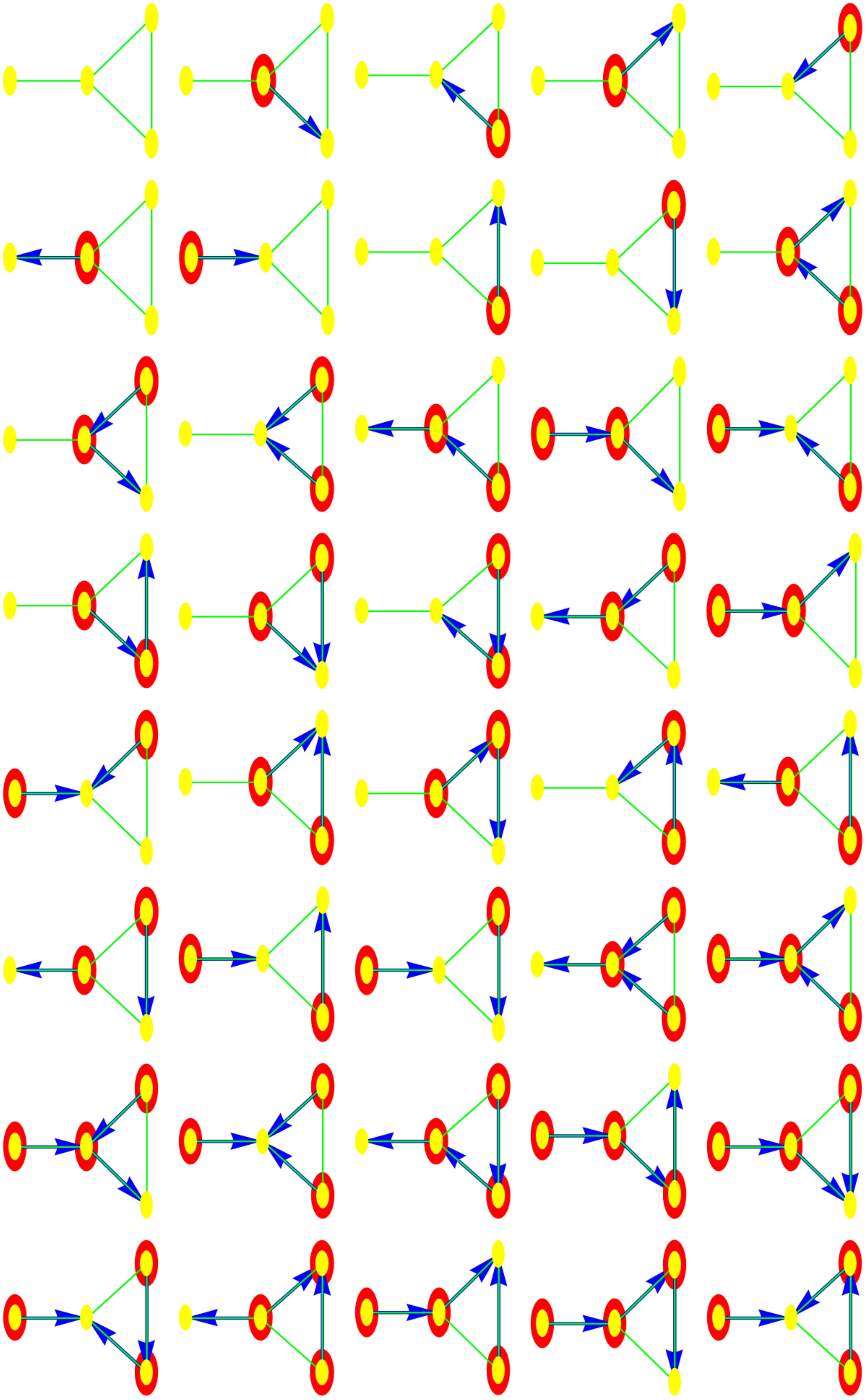}}
\end{figure}

{\bf Example 3.}
The kite graph $G$ has $L=\left[
                 \begin{array}{cccc}
                  3 & -1 & -1 & -1 \\
                  -1 & 2 & -1 & 0 \\
                  -1 & -1 & 3 & -1 \\
                  -1 & 0 & -1 & 2 \\
                 \end{array}
                 \right]$ with eigenvalues $4,4,2,0$
so that $f(G) =75$. The last $\Det(L)=32$ forests match rooted maximal spanning trees.

\begin{figure}[H]
\scalebox{0.26}{\includegraphics{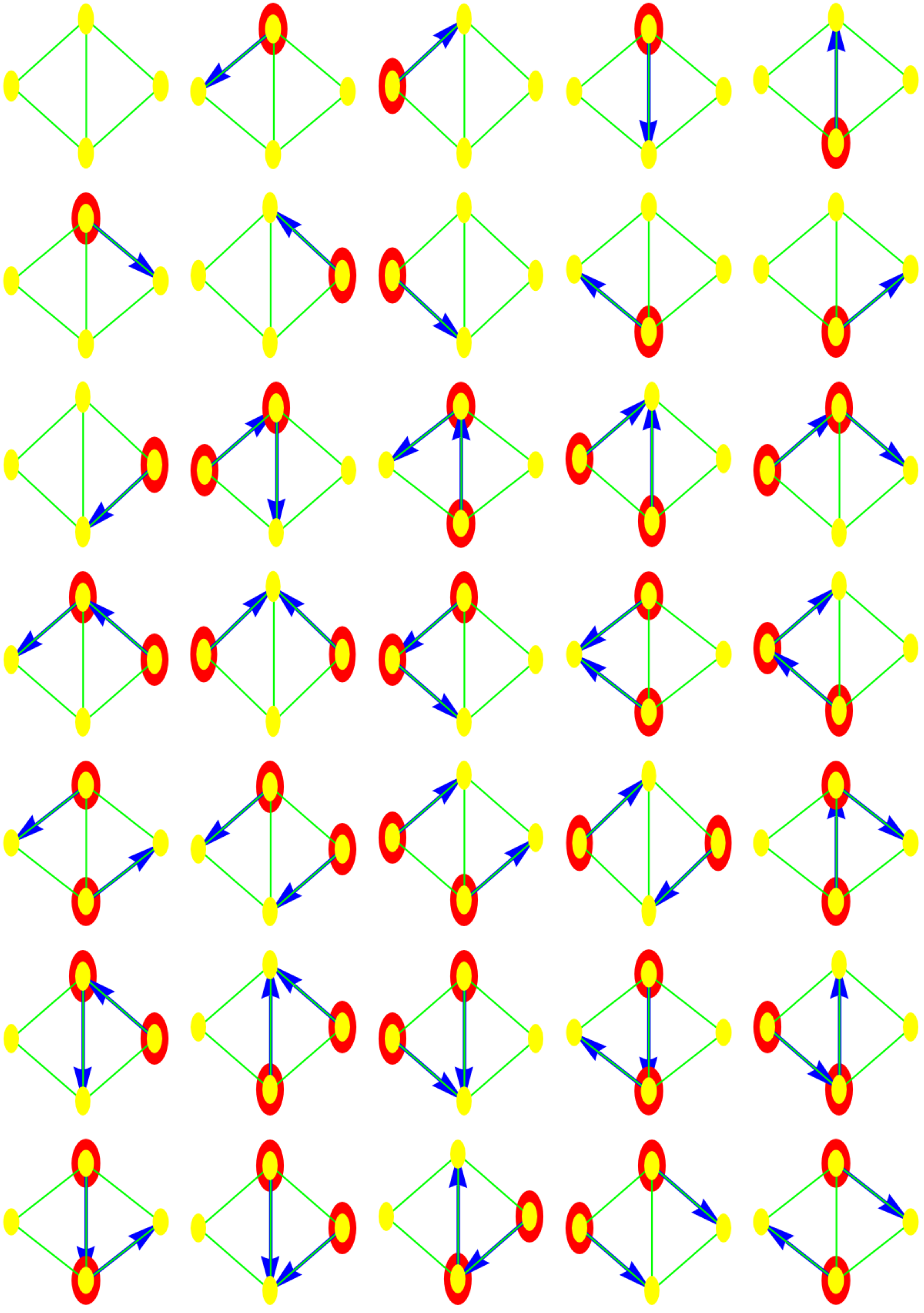}}
\end{figure}
\begin{figure}
\scalebox{0.32}{\includegraphics{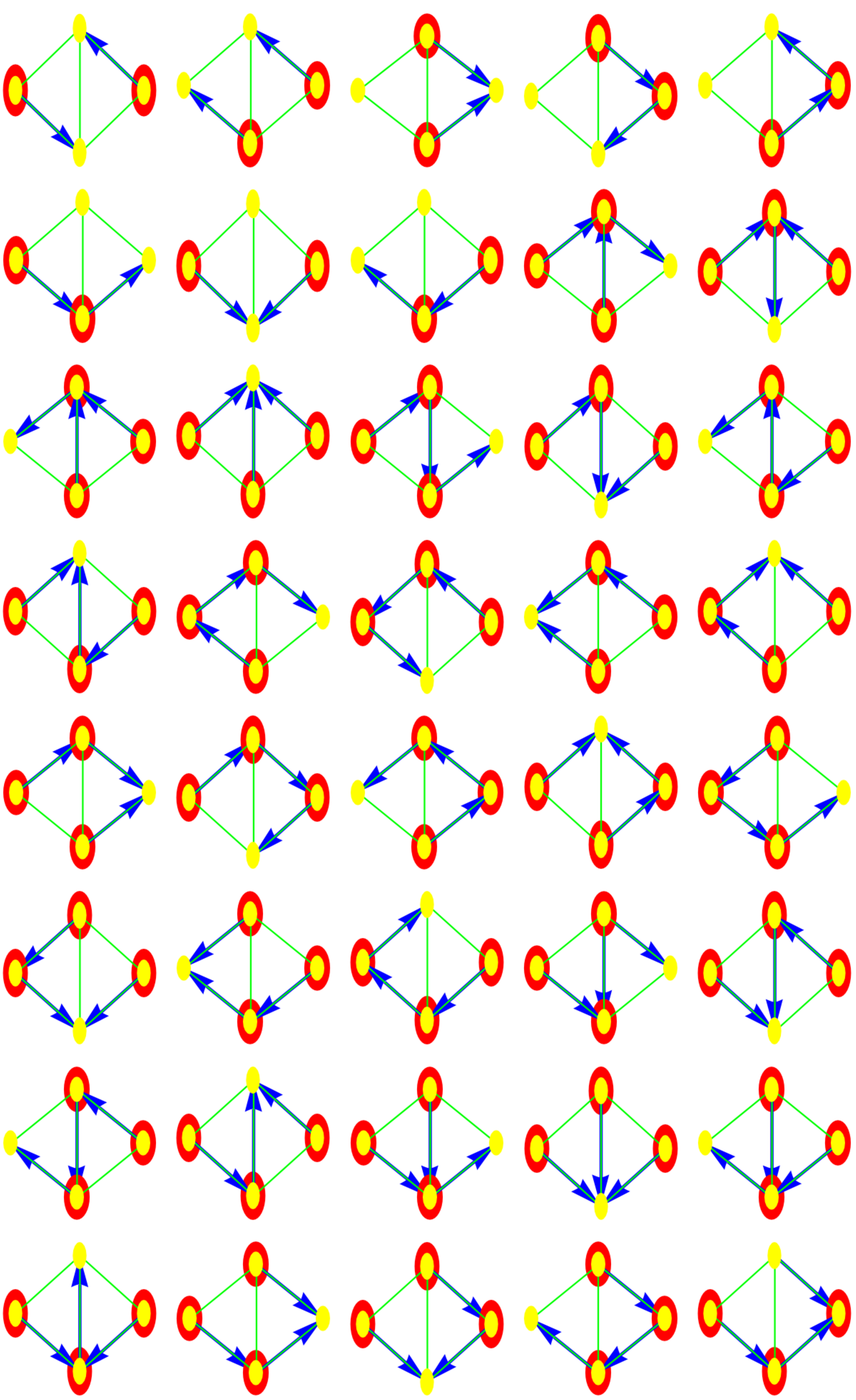}}
\end{figure}

{\bf Example 4.}
The tadpole graph $G$ has $L= \left[
                 \begin{array}{ccccc}
                  2 & -1 & 0 & -1 & 0 \\
                  -1 & 2 & -1 & 0 & 0 \\
                  0 & -1 & 2 & -1 & 0 \\
                  -1 & 0 & -1 & 3 & -1 \\
                  0 & 0 & 0 & -1 & 1 \\
                 \end{array}
                 \right]$ with eigenvalues $4.4812.., 2.6889.., 2, 0.8299..,0$
and $f(G) =111$. The last $\Det(L)=20$ match rooted spanning trees.

\begin{figure}[H]
\scalebox{0.222}{\includegraphics{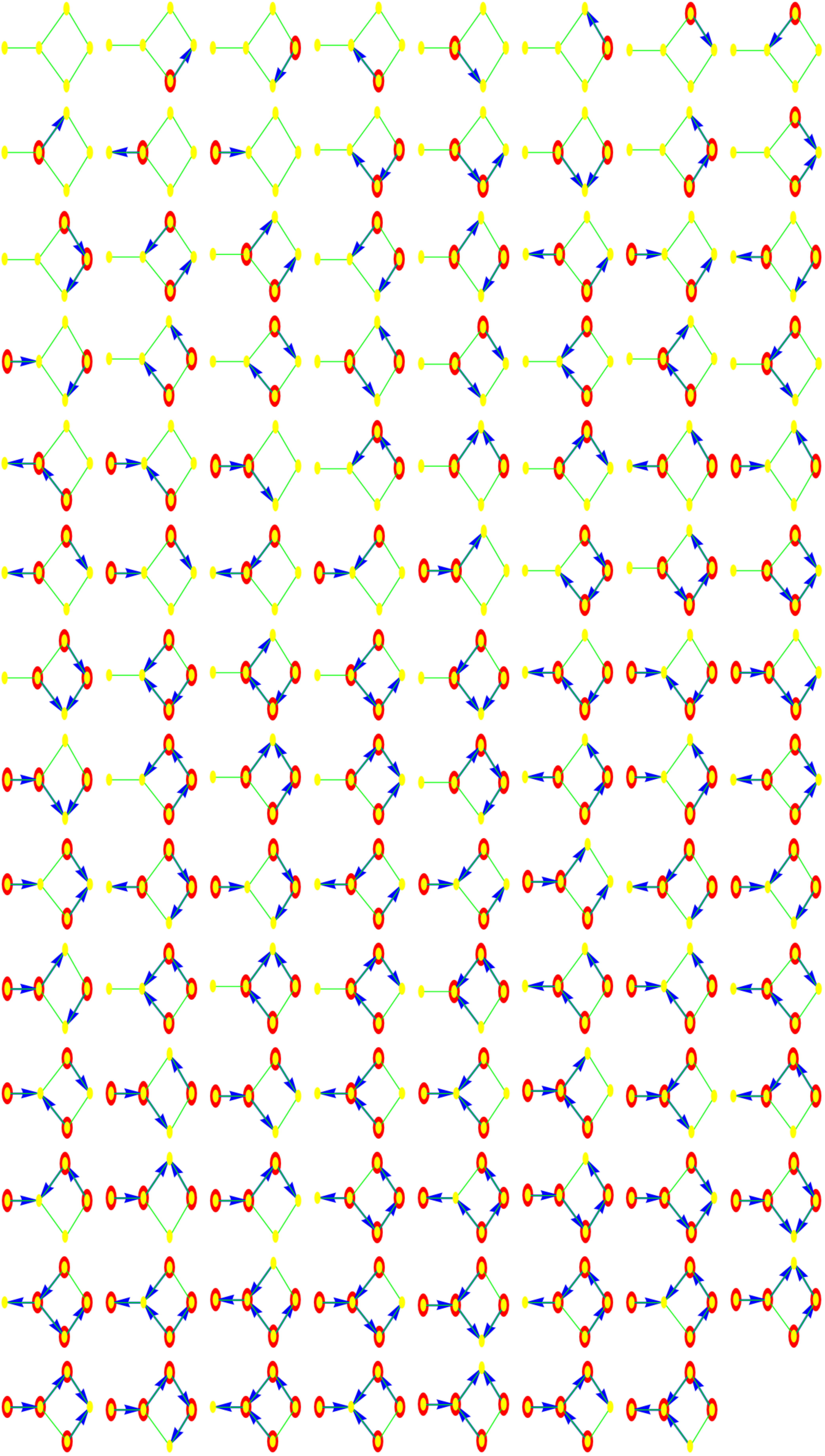}}
\end{figure}

{\bf Example 5.}
The extended complete graph $G=K_4^+$ has 
\begin{tiny}$L= \left[
                 \begin{array}{ccccc}
                  3 & -1 & -1 & -1 & 0 \\
                  -1 & 3 & -1 & -1 & 0 \\
                  -1 & -1 & 3 & -1 & 0 \\
                  -1 & -1 & -1 & 4 & -1 \\
                  0 & 0 & 0 & -1 & 1 \\
                 \end{array}
                 \right]$ \end{tiny} with eigenvalues $5,4,4,1,0$
and $f(G) =300$. The last $\Det(L)=80$ match rooted spanning trees.

\begin{figure}[H]
\scalebox{0.28}{\includegraphics{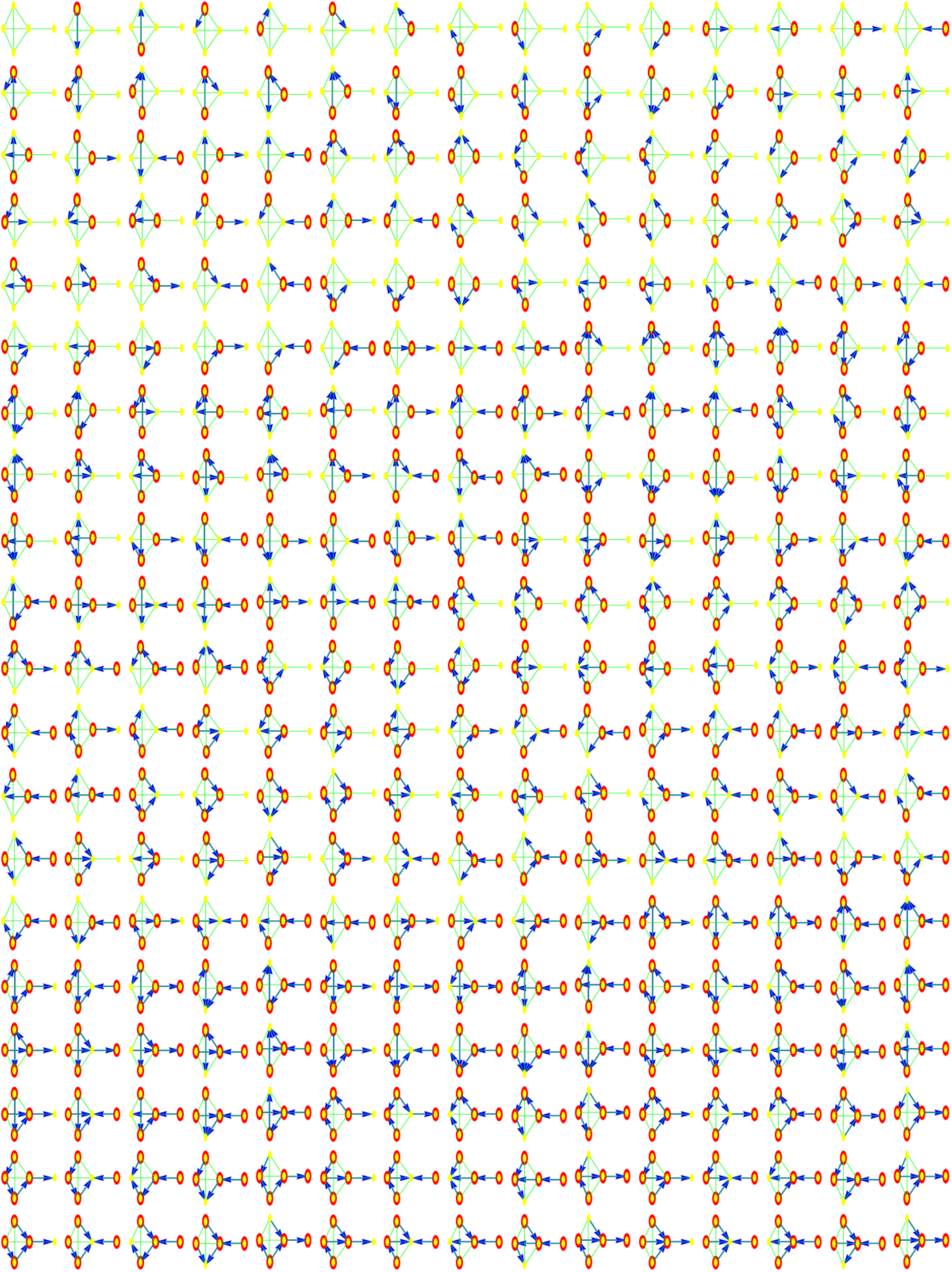}}
\end{figure}

{\bf Example 6.}
The graph $G$ with Laplacian
\begin{tiny}$L= \left[
 \begin{array}{cccccc}
                  4 & -1 & -1 & -1 & 0 & -1 \\
                  -1 & 3 & -1 & -1 & 0 & 0 \\
                  -1 & -1 & 3 & -1 & 0 & 0 \\
                  -1 & -1 & -1 & 4 & -1 & 0 \\
                  0 & 0 & 0 & -1 & 2 & -1 \\
                  -1 & 0 & 0 & 0 & -1 & 2 \\
                 \end{array}
                 \right]$ \end{tiny} has the spectrum $\{4+\sqrt{2},3+\sqrt{3},4,4-\sqrt{2},3-\sqrt{3},0 \}$
and $f(G) =1495$. The last $\Det(L)=336$ forests match rooted spanning trees.

\begin{figure}[H]
\scalebox{0.24}{\includegraphics{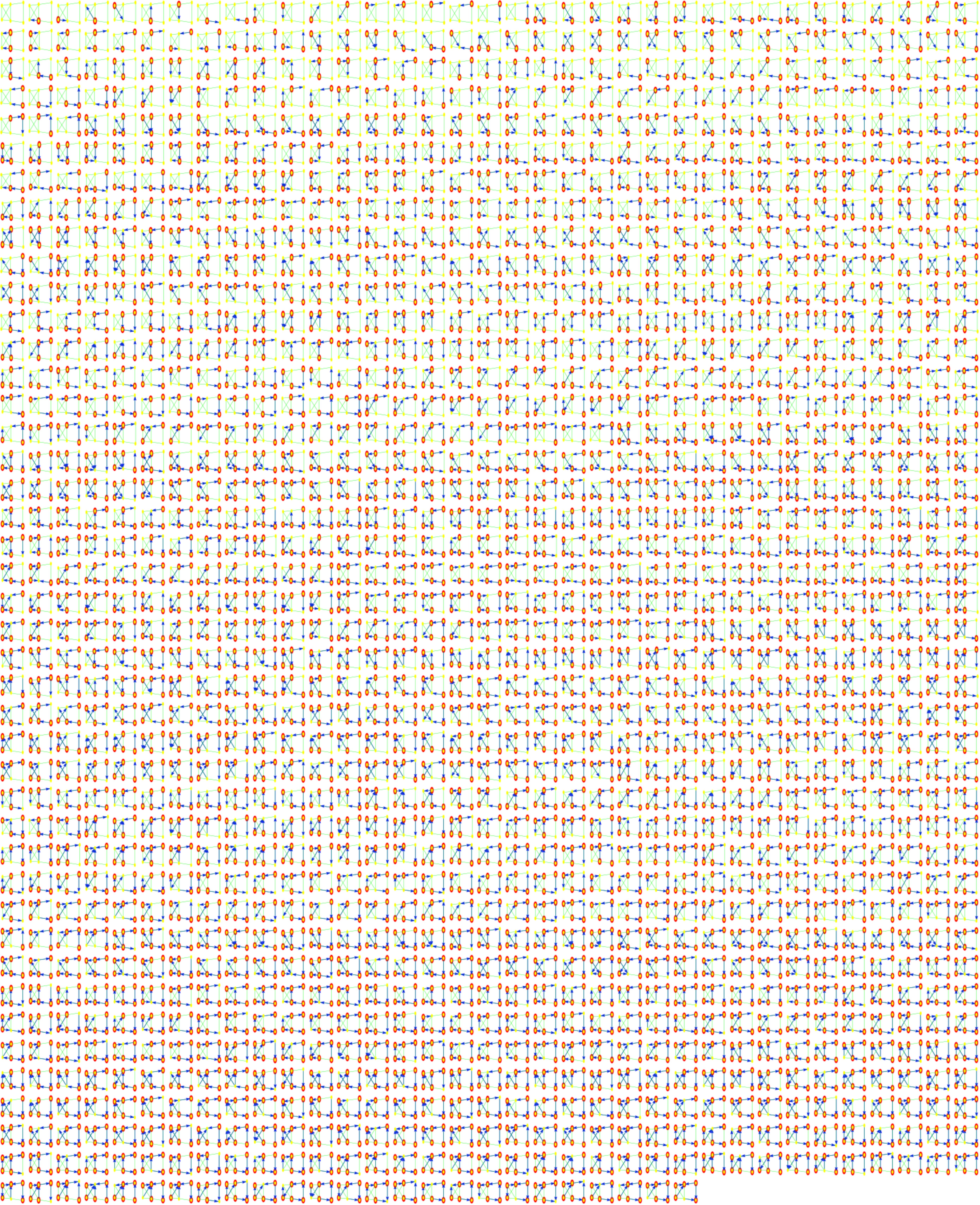}}
\end{figure}

\pagebreak 

\end{document}